\documentclass[9pt]{article}
\usepackage{latexsym,amsfonts,euscript,amsthm}

\usepackage{amssymb}
\usepackage{stmaryrd}
\usepackage{mathrsfs,amsmath}
\usepackage{leftidx}

\usepackage[figuresright]{rotating}
\usepackage[table]{xcolor}
\usepackage{multirow}
\usepackage{booktabs}
\usepackage{tabularx}

\usepackage{tocbibind} 

\usepackage[pagebackref]{hyperref}

\hypersetup{%
  colorlinks=true,
linktoc=all,
  citebordercolor=green,
  linkbordercolor=red,
 citecolor=green,
    linkcolor=red,
  urlcolor=blue
}

\newtheorem{lem}{\bf Lemma}[section]

\newtheorem{thm}[lem]{\bf Theorem}
\newtheorem{rmk}[lem]{\bf Remark}

\newtheorem{mainthm}{Theorem}

\newcommand{\ackname}{Acknowledgements}
\makeatletter
\if@titlepage
  \newenvironment{acknowledgement}{%
    \titlepage
    \null\vfil
    \@beginparpenalty\@lowpenalty
    \begin{center}%
      \bfseries \ackname
      \@endparpenalty\@M
    \end{center}}%
  {\par\vfil\null\endtitlepage}
\else
  \newenvironment{acknowledgement}{%
    \if@twocolumn
      \section*{\ackname}%
    \else
      \small
      \begin{center}%
        {\bfseries \ackname\vspace{-0.5em}\vspace{\z@}}%
      \end{center}%
      \quotation
    \fi}
    {\if@twocolumn\else\endquotation\fi}
\fi
\makeatother

\makeatletter 
\@addtoreset{equation}{section}
\makeatother  

\makeatletter
\def\thanks#1{\protected@xdef\@thanks{\@thanks\protect\footnotetext{#1}}}

				   \AtEndDocument{\bigskip{\footnotesize
				   \textsc{Department of Mathematics, Suzhou University of Technology, No.99 South Third Ring
				   Road,
				   Changshu, Jiangsu, 215500, China.} \par  
				   \textit{Email address}: \href{mailto:yuzeng2004@163.com}{yuzeng2004@163.com} \par
				   \addvspace{\medskipamount}
				   \textsc{School of Mathematics and Physics, Suqian University, No. 399 South Yellow River Road,
				   Suqian, Jiangsu, 223800, China.} \par  
				   \textit{Email address}: \href{mailto:leejinbao25@163.com,leejinbao25@hotmail.com}{leejinbao25@163.com}   \par
				   \addvspace{\medskipamount}
				    \textsc{Department of Mathematics, Texas State University, 601 University Drive, San Marcos, TX 78666, USA.} \par  
				   \textit{Email address}: \href{mailto:yang@txstate.edu}{yang@txstate.edu} \par
				 }}

\title{On $p$-parts of conjugacy class sizes and the index of the $p$-core
\thanks{
\textbf{2020 MR Subject Classification}\,\, Primary 20E45, 20D25, Secondary 20C15\\
The first author is supported by the  NSF of China (No.12171058, 12301018), the Natural Science Foundation of Jiangsu Province (No. BK20231356) and the Natural Science Foundation of the Jiangsu Province Higher Education Institutions of China (No. 23KJB110002).
The second author is supported by the NSF of China (No. 12471017), the Scientific Research Foundation for Advanced Talents of Suqian University (No. 2022XRC069), the Natural Science Foundation of Suqian, China (No. K202432).
The third author was partially supported by a grant from the Simons Foundation (No. 918096).
}}

\author{Yu Zeng, Jinbao Li*\thanks{*Corresponding author}, Yong Yang\\
\small{\emph{Dedictaed to Professor Guohua Qian on his 60th birthday}}}

\date{}

\begin{document}
\maketitle

\begin{abstract}
   For a prime $p$ and an arbitrary finite group $G$,
   we show that if $p^{2}$ does not divide the size of each conjugacy class of \emph{$p$-regular} element (element of order not divisible by $p$)
   in $G$, then the largest power of $p$ dividing the index $|G:\mathbf{O}_{p}(G)|$ is at most $p^{2}$.
\end{abstract}

\section{Introduction}

 The fundamental It\^o-Michler theorem states that, for a finite group $G$, the prime $p$ does not divide the degree of every irreducible complex character of $G$ if and only if $G$ has an abelian normal Sylow $p$-subgroup.
 Specifically, in this scenario, the index $|G:\mathbf{O}_{p}(G)|$ is not divisible
 by $p$, where $\mathbf{O}_{p}(G)$ denotes the \emph{$p$-core} of $G$ (the unique largest normal $p$-subgroup of $G$).

Lewis, Navarro et al. \cite{lewis15,lewis14} introduced a direction of generalizing the It\^o-Michler theorem.
Suppose that the $a$-th power of a prime $p$ does not divide the degree of every irreducible complex character (resp. irreducible $p$-Brauer character) of a finite group $G$.
Is it possible to bound $|G:\mathbf{O}_{p}(G)|_p$
(the largest power of $p$ dividing $|G:\mathbf{O}_{p}(G)|$)? 
This question has been investigated by several authors, including \cite{bessenrodt20,lewis15,lewis14,li21,qian18,qian25,yang18b}.
For instance, Michler \cite{michler86} classified such finite groups $G$ when $a=1$;
recently, Qian \cite{qian25} established the sharp upper bound for $|G:\mathbf{O}_{p}(G)|_p$ when $a=2$.

A dual version of this question arises in the context of conjugacy class sizes. 
Suppose that the $a$-th power of a prime $p$ does not divide the size of the conjugacy class of every $p$-regular element in a finite group $G$.
Is it possible to bound $|G:\mathbf{O}_{p}(G)|_p$?  
This question has been studied by several authors such as \cite{bessenrodt20,camina72,li19,yang18a,yang18b}.
Camina \cite{camina72} classified such finite groups $G$ when $a=1$.
When $a=2$,
Yang and Qian \cite{yang18a} showed that $|G:\mathbf{O}_{p}(G)|_p\leq p^4$ and also provided the sharp upper bound for $|G:\mathbf{O}_{p}(G)|_p$ when $p=2$;
Li and Yang \cite{li19} subsequently improved the upper bound for $|G:\mathbf{O}_{p}(G)|_p$ to $p^{3}$ for odd primes $p$.
Inspired by a recent result of Qian \cite[Theorem 1.1]{qian25}, we establish in the following theorem that when $a=2$, the sharp upper bound for $|G:\mathbf{O}_{p}(G)|_p$ is $p^a$.

\begin{mainthm}\label{thmA}
	Let $G$ be a finite group.
	If $p^{2}$ does not divide $|G:\mathbf{C}_{G}(g)|$ for each $p$-regular element $g$ of $G$,
	then $|G:\mathbf{O}_{p}(G)|_p\leq p^{2}$. 
  \end{mainthm}

\begin{rmk}
 {\rm  
 The example in \cite[Example 3.3]{lewis14} demonstrates that our bound in Theorem \ref{thmA} is sharp.
 For readers' convenience, we present the example here.
 
 By Dirichlet's theorem \cite[Page 169, (b)]{isaacs76}, there exists a prime $q$ such that $p^2$ divides $q-1$. 
 Note that, by \cite[Lemma 10.11]{isaacs76}, $ (q^p-1 ) /(q-1)$ is divisible by $p$, but not by $p^2$. 
 Let $\Gamma=\Gamma (q^p )$ denote
 the semilinear group acting on its natural module $V$ of order $q^p$. 
 Then $\Gamma$ contains a normal cyclic subgroup $\Gamma_0$ of order $q^p-1$, which acts fixed-point-freely on $V$.
 Additionally, $\Gamma / \Gamma_0$ has order $p$, and $\Gamma$ splits over $\Gamma_0$. Let $x \in \Gamma$ have order $p$ not contained in $\Gamma_0$. 
 We may assume the automorphism induced by $x$ on $\Gamma_0$ is $a \mapsto a^q$. 
 The centralizer of $x$ in $\Gamma_0$ (i.e. $\mathbf{Z}(\Gamma)$) has order $q-1$. 
 Now, $a \mapsto[a, x]$ maps $\Gamma_0$ onto $U:= [\Gamma_0, x ]$. 
 In particular, $ |\Gamma_0: U |=q-1$ and $\Gamma / U=\Gamma_0 / U \times H / U$ where $H=\langle U, x\rangle=U \rtimes \langle x\rangle$ has order $p(q^{p}-1)/(q-1)$. 
 So, $H$ is normal in $\Gamma$, implying $H$ contains every $\Gamma$-conjugate of $x$. 
 For a nontrivial $v \in V$, we have that $\mathbf{C}_{\Gamma}(v)= \langle x^y \rangle$ for some $y \in \Gamma$,
  so $\mathbf{C}_H(v)$ has order $p$.

  Set $G=V \rtimes H$.
  Then $\mathbf{O}_{p}(G)=1$ and $|G|_p=p^{2}$.
  Let $g$ be a nontrivial $p$-regular element of $G$.
  As $V\mathbf{O}_{p'}(H)$ is a normal Hall $p'$-subgroup of the solvable group $G$, 
  it forces that $g\in V \mathbf{O}_{p'}(H)$.
  Note that $V\mathbf{O}_{p'}(H)$ is a Frobenius group with kernel $V$,
  and so either $g\in V$, or $g^v\in \mathbf{O}_{p'}(H)$ for some $v\in V$. 
  If the former holds,
  as $\mathbf{C}_{G}(g)=V\mathbf{C}_{H}(g)$, we deduce that $|\mathbf{C}_{G}(g)|_p=|\mathbf{C}_{H}(g)|=p$.
  If the latter holds,
  observing that $[\mathbf{O}_{p}(H),\mathbf{O}_{p'}(H)]=1$ and that $|\mathbf{O}_{p}(H)|\geq p$, we conclude that $|\mathbf{C}_{G}(g)|=|\mathbf{C}_{G}(g^v)|\geq p$.
  Consequently, $p^{2}\nmid |G:\mathbf{C}_{G}(g)|$ for every $p$-regular element $g$ of $G$.

  Note that one can also construct a nonsolvable example if $p>2$.
  For instance, the direct product of $G$ (defined in the last paragraph) and a nonabelian simple $p'$-group.
 } 
\end{rmk}

Throughout the rest of this paper, we only consider finite groups and complex characters.

\section{Proofs}

We use
\cite{huppert67} as a source for results in group theory and
\cite{isaacs76} for results in character theory.
Throughout the paper, we consistently refer $p$ as a prime.
For a positive integer $n$ and a prime $p$, we write $n_p$ to
denote the largest $p$-power divisor of $n$.
 Other notation will be recalled or defined when necessary.

We begin with an elementary result concerning conjugacy classes.

\begin{lem}\label{lem: conj}
	Let $N$ be a normal subgroup of a finite group $G$, and let $\overline{G}=G/N$.
  \begin{description}
	\item[(1)] If $x \in N$, then $|x^N|$ divides $|x^G|$.
	\item[(2)] If $x \in G$, then $|\overline{x}^{\overline{G}}|$ divides $|x^G|$.
  \end{description}
\end{lem}
\begin{proof}
	Since $N$ is normal in $G$, $\mathbf{C}_{G}(g) N$ is a subgroup of $G$ for every $g\in G$.
    So, for $x\in N$, 
	$$|x^{N}|=|N:\mathbf{C}_{N}(x)|=|\mathbf{C}_{G}(x)N:\mathbf{C}_{G}(x)| ~\big|~ |G:\mathbf{C}_{G}(x)|=|x^{G}|.$$
    Note that $\overline{\mathbf{C}_{G}(g)}\leq \mathbf{C}_{\overline{G}}(\overline{g})$ for every $g\in G$,
	and so, for $x\in G$,
   \[
   |\overline{x}^{\overline{G}}|=|\overline{G}:\mathbf{C}_{\overline{G}}(\overline{x})| ~\big|~ |G:\mathbf{C}_{G}(x)N|~\big|~  |G:\mathbf{C}_{G}(x)|=|x^{G}|. 
   \]
\end{proof}

Let $G$ be a finite group, $p$ a prime and $n$ a positive integer.
For the sake of brevity, we say that $G$ satisfies $\mathbf{SCpR}(p^{n})$ (or that
$G\in \mathbf{SCpR}(p^{n})$, for short) if $p^{n}$ does not divide $|G:\mathbf{C}_{G}(g)|$ for each $p$-regular 
element $g$ of $G$.
The next lemma shows that the property $\mathbf{SCpR}(p^{n})$ is inherited by normal subgroups and factor groups.

\begin{lem}\label{lem: ind}
  Let $G\in \mathbf{SCpR}(p^{n})$ and $N$ a normal subgroup of $G$.
 Then both $N$ and $G/N$ satisfy $\mathbf{SCpR}(p^{n})$.
\end{lem}
\begin{proof}
	Let $x$ be a $p$-regular element in $N$.
	Since $|x^N|$ divides $|x^G|$ by Lemma \ref{lem: conj}, we deduce that $N\in \mathbf{SCpR}(p^{n})$.
	Set $\overline{G}=G/N$.
	Let $\overline{x}$ be a $p$-regular element in $\overline{G}$.
	Write $o(x)=p^a m$ where $(p,m)=1$, and set $y=x^{p^a}$.
	Then $\langle \overline{x}\rangle=\langle \overline{y}\rangle$, implying that $\mathbf{C}_{\overline{G}}(\overline{x})=\mathbf{C}_{\overline{G}}(\overline{y})$.
	Noting that $|\overline{x}^{\overline{G}}|=|\overline{y}^{\overline{G}}|$ divides $|y^G|$ by Lemma \ref{lem: conj},
	we conclude that $\overline{G}\in \mathbf{SCpR}(p^{n})$.
\end{proof}

\begin{thm}\label{thm: qian}
	Let $G$ be a finite group and $p$ an odd prime.
	If $p^{2}$ does not divide $\chi(1)$ for each $\chi \in \mathrm{Irr}(G)$,
	then $|G:\mathbf{O}_{p}(G)|_p\leq p^{2}$. 
\end{thm}
\begin{proof}
	This is a partial result of \cite[Theorem 1.1]{qian25}.
\end{proof}

Recall that the solvable radical of a finite group $G$ is the unique largest normal 
solvable subgroup of $G$.

\begin{lem}\label{lem: sol=1}
	 Let $G\in \mathbf{SCpR}(p^2)$ for an odd prime $p$.
	 If $G$ has a trivial solvable radical, then $|G|_p\leq p$.
\end{lem}
\begin{proof}
  This is \cite[Lemma 2.9]{li19}.
\end{proof}

Let $p$ be a prime, $B$ a $p$-block of a finite group $G$ with a defect group $D$.
Then $\mathrm{Irr}(G)$ (the set of irreducible characters of $G$) is a disjoint union of $\mathrm{Irr}(B)$ (the \emph{set of irreducible characters in $B$}) with $B$ running over all $p$-blocks of $G$.
For $\chi\in \mathrm{Irr}(B)$, it is well-known that 
\[
 \chi(1)_p=p^{\mathrm{ht}(\chi)}|G:D|_p,
\]
 where $\mathrm{ht}(\chi)$ denotes the \emph{height} of $\chi$.

\begin{lem}[\cite{brauer41}]\label{lem: height of character}
	Let $B$ be a $p$-block of a finite group $G$ with a nonabelian defect group $D$.
	If $\chi \in \mathrm{Irr}(B)$,
	then $p^{\mathrm{ht}(\chi)}\leq |D|/p^{2}$.
\end{lem}

Now, we are ready to prove Theorem \ref{thmA}.

\begin{proof}[Proof of Theorem \ref{thmA}]
	  Let $G$ be a counterexample of minimal possible order.
  Then $|G:\mathbf{O}_{p}(G)|_p>p^{2}$.
  Applying \cite[Theorem 1.2]{yang18a}, we deduce that $p>2$.
  Since $G/\mathbf{O}_{p}(G)\in \mathbf{SCpR}(p^{2})$ by Lemma \ref{lem: ind}, by the minimality of $G$,
 the assumption that $\mathbf{O}_{p}(G)>1$ would lead to a contradiction that $|G:\mathbf{O}_{p}(G)|_p\leq p^{2}$.
 So, we may assume that $\mathbf{O}_{p}(G)=1$.
 Let $K$ be the solvable radical of $G$.
Since $G/K\in \mathbf{SCpR}(p^2)$ has a trivial solvable radical,
  Lemma \ref{lem: sol=1} implies that $|G/K|_p\leq p$.
  Note that $K\in \mathbf{SCpR}(p^2)$ is the solvable radical of $G$, and hence \cite[Theorem 1.1]{yang18a} yields that $|K|_p=|K:\mathbf{O}_{p}(G)|_p\leq p^{2}$.
  As $|G|_p>p^{2}$, it forces $|G|_p=p^{3}$.

  Let $\chi\in \mathrm{Irr}(G)$.
	Then there exists a $p$-block $B$ of $G$ with defect group $D$ such that $\chi \in \mathrm{Irr}(B)$.
	Note that $D$ is a Sylow $p$-subgroup of $\mathbf{C}_{G}(x)$ for some $p$-regular element $x \in G$,
	and so $G\in \mathbf{SCpR}(p^2)$ implies that $|G:D|_p\leq p$.
	As $|G|_p= p^{3}$, it forces $|D|=p^{2}$ or $p^{3}$.
	If $D$ is abelian, as by \cite[Theorem 1.1]{kessar13} $\mathrm{ht}(\chi)=0$, we deduce that $\chi(1)_p=p^{\mathrm{ht}(\chi)}|G:D|_p \leq p$.
    If $D$ is nonabelian, then $|D|=|G|_p=p^{3}$.
	As $p^{\mathrm{ht}(\chi)}\leq |D|/p^{2}$ by Lemma \ref{lem: height of character}, 
	we deduce that $\chi(1)_p=p^{\mathrm{ht}(\chi)}|G:D|_p \leq p$.

	Applying Theorem \ref{thm: qian} leads to the contradiction that $|G|_p\leq p^{2}$.
\end{proof}

\begin{acknowledgement}
	The authors are grateful to the referee for her/his
 valuable comments.
\end{acknowledgement}


\begin{thebibliography}{999}\setlength{\itemsep}{-2mm} 
\small
	\bibitem{bessenrodt20}
	C. Bessenrodt and Y. Yang,
	\newblock On $p$-parts of Brauer character degrees and $p$-regular conjugacy class sizes of finite groups,
	\newblock \emph{J. Algebra} {\bf 560} (2020), 296--311.


	\bibitem{brauer41}
	R. Brauer,
	\newblock Investigations on group characters,
	\newblock \emph{Ann. of Math.} {\bf 42(2)} (1941), 936--958.



	\bibitem{camina72}
	A.R. Camina,
	\newblock Arithmetical conditions on the conjugacy class numbers of a finite
	group,
	\newblock \emph{J. Lond. Math. Soc.} {\bf 5(2)} (1972), 127--132.



	\bibitem{huppert67}
	B. Huppert,
	\emph{Endliche gruppen I}, Springer Verlag, Berlin, 1967.


\bibitem{kessar13}
	R. Kessar and G. Malle,
	\newblock Quasi-isolated blocks and Brauer's height zero conjecture,
	\newblock \emph{Ann. of Math.} {\bf 178(2)} (2013), 321--384.


	\bibitem{isaacs76}
	I.M. Isaacs,
	\emph{Character Theory of Finite Groups}, Academic Press, New York, 1976.







	\bibitem{lewis15}
	M.L. Lewis, G. Navarro, P.H. Tiep and H.P. Tong-Viet,
	\newblock $p$-parts of character degrees,
	\newblock \emph{J. Lond. Math. Soc.} {\bf 92(2)} (2015), 483--497.

	\bibitem{lewis14}
	M.L. Lewis, G. Navarro and T.R. Wolf,
	\newblock $p$-parts of character degrees and the index of
	the Fitting subgroup,
	\newblock \emph{J. Algebra} {\bf 411} (2014), 182--190.





	\bibitem{li19}
	J. Li and Y. Yang,
	\newblock A note on $p$-parts of conjugacy class sizes,
	\newblock \emph{J. Group Theory} {\bf 22(5)} (2019), 933--940.

	
	\bibitem{li21}
	J. Li and Y. Yang,
	\newblock A note on $p$-parts of Brauer character degrees,
	\newblock \emph{Bull. Aust. Math. Soc.} {\bf 103(1)} (2021), 83--87.



	\bibitem{michler86}
	G.O. Michler,
	\emph{Brauer's conjectures and the classification of finite simple
 groups}, Lect. Notes Math. 1178, Springer, Berlin, 1986.




	
		  



	\bibitem{qian18}
	G. Qian,
	\newblock A note on $p$-parts of character degrees,
	\newblock \emph{Bull. Lond. Math. Soc.} {\bf 50(4)} (2018), 663--666.

	\bibitem{qian25}
	G. Qian,
	\newblock On $p$-parts of character degrees,
	\newblock \emph{J. Pure Appl. Algebra} {\bf 229(1)} (2025), 107793.






	\bibitem{yang18a}
	Y. Yang and G. Qian,
	\newblock On $p$-parts of conjugacy class sizes of finite groups,
	\newblock \emph{Bull. Aust. Math. Soc.} {\bf 97(3)} (2018), 406--411.

    \bibitem{yang18b}
	Y. Yang and G. Qian,
	\newblock On $p$-parts of character degrees and conjugacy class sizes of finite groups,
	\newblock \emph{Adv. Math.} {\bf 328} (2018), 356--366.




\end{thebibliography}
\end{document}